\documentclass[12pt,a4paper,reqno]
{amsart}
\usepackage{amsfonts, amsmath, amssymb, amsthm, amsbsy, amscd, euscript}

\usepackage[english]{babel}

\usepackage[top=2.25cm, bottom=2.25cm, left=2.25cm, right=2.25cm]{geometry}

\usepackage{stmaryrd}

\usepackage{tikz,tikz-cd}
\usetikzlibrary{arrows,positioning,automata,shadows,fit,shapes}

\usepackage{cmap}

\usepackage[pdftex,
unicode=true,
colorlinks=false
]{hyperref}

\usepackage{DLdef1}

\renewcommand {\ssbegin}[2][*]
 {\refstepcounter{subsection}%
\if#1*
\addcontentsline{toc}{subsection}{\thesubsection.\hskip 1pc #2}%
\else
\addcontentsline{toc}{subsection}{\thesubsection.\hskip 1pc #2. #1}%
\fi
 \def \secno {\gdef \secno {}{\ssecfont
\thesubsection.\hskip 2ex}%
 }%
 \begin{#2}}

\renewcommand {\sssbegin}[2][*]
  {\refstepcounter{subsubsection}
\if#1*
\addcontentsline{toc}{subsubsection}{\thesubsubsection.\hskip 1pc #2}%
\else
\addcontentsline{toc}{subsubsection}{\thesubsubsection.\hskip 1pc #2. #1}
\fi
  \def \secno {\gdef \secno {}{\ssecfont \thesubsubsection.\hskip 2ex}%
  }%
   \begin{#2}}

\renewcommand {\parbegin}[2][*]
  {\refstepcounter{paragraph}
\if#1*
\addcontentsline{toc}{paragraph}{\theparagraph.\hskip 1pc #2}%
\else
\addcontentsline{toc}{paragraph}{\theparagraph.\hskip 1pc #2. #1}
\fi
  \def \secno {\gdef \secno {}{\ssecfont \theparagraph.\hskip 2ex}%
  }%
   \begin{#2}}

\renewcommand{\ssecfont}{\normalfont}

\newcommand{\dd}{\mathrm{d}}

\newcommand{\del}{\partial}

\DeclareMathOperator{\Cl}{Cl}

\newcommand{\twedge}{{\textstyle{\bigwedge}}}

\newcommand{\fX}{\mathfrak{X}}

\newcommand{\poly}{\mathrm{poly}}
\newcommand{\dyn}{\mathrm{dyn}}

\title{Group valued moment maps for even and odd simple $G$-modules}

\author{Anton Alekseev}
\address[A.~Alekseev]{Section of Mathematics, University of Geneva, Rue du Conseil G\'en\'eral 7-9, 12211, Geneva, Swtizerland} 
\email{Anton.Alekseev@unige.ch}

\author{Andrey Krutov}
\address[A.~Krutov]{Mathematical Institute of Charles University, Sokolovsk\'a 83, 186 75 Prague, Czech Republic} 
\email{andrey.krutov@matfyz.cuni.cz}



\begin{document}

\begin{abstract}
  Let~$G$ be a~complex simple Lie group, and~$\mathfrak{g}$ its Lie algebra. It is well known that
  a finite-dimensional $G$-module~$V$  carrying  a~nondegenerate invariant bilinear form gives rise to
  a~Hamiltonian Poisson space with a~quadratic moment map~$\mu$.
  We show that under condition $\mathrm{Hom}_\fg({\textstyle{\bigwedge}}^3 V, S^3V)=0$ this space can be viewed as a~quasi-Poisson space with the same bivector, and with the group valued moment map $\Phi = \exp \circ \mu$. Furthermore, we show that by modifying the bivector by the standard $r$-matrix for $\mathfrak{g}$ one obtains a space with a Poisson action of the Poisson-Lie group~$G$, and with the moment map in the sense of Lu taking values in the dual Poisson-Lie group~$G^\ast$.
\end{abstract}

\maketitle

\section{Introduction}

Let $G$ be a~complex simple Lie group, $\fg$ be its Lie algebra, and $V$ be a~finite-dimensional $G$-module equipped with an
invariant holomorphic symplectic form.
Then, $V$ becomes a Hamiltonian Poisson space with Poisson structure $\pi$ given by the inverse of the symplectic form, and with the quadratic moment map $\mu$ defined by the $\fg$-action on $V$. 

Similarly, if $V$ carries a $G$-invariant nondegenerate symmetric bilinear form, then following
Kostant--Sternberg~\cite{KostantSternberg} one can define $G$-invariant symplectic and Poisson structures
(again denoted by $\pi$) on the supermanifold $\Pi V$, where $\Pi$ stands for the parity change functor. As before, the $G$-action is Hamiltonian, and the moment map $\mu$ is quadratic.

In this paper, we restrict our attention to $G$-modules with a rather restrictive extra property
\begin{equation}      \label{eq:key}
  \Hom_\fg(\twedge^3 V, S^3(V))=0.
\end{equation}
Such modules can be classified (see Section 3 for details), and by the results of \cite{Zwicknagl2009} they carry quadratic $r$-matrix Poisson brackets $r_V$. Furthermore, the  Poisson-Lie group $G$ defined by the $r$-matrix acts on
$(V, r_V)$ (or $(\Pi V,  r_{\Pi V})$).

Our main results can be summarised in the following two theorems: 
\begin{Theorem}
  Let $V$ be a~finite-dimensional $G$-module which satisfies condition~\eqref{eq:key}, and which carries a~$G$-invariant symplectic form.
  Then, the corresponding Poisson space $(V, \pi)$ can be viewed as a~Hamiltonian quasi-Poisson space with the
  same bivector, and with the group valued moment map $\Phi = \exp \circ \mu: V \to G$.
  Furthermore, the space $(V, \pi-r_{V})$ is Poisson, and it carries an action of the Poisson-Lie group~$G$
  equipped with the standard Poisson bracket.
  The pre-image under $\Phi$ of the domain of the Gauss decomposition carries a Hamiltonian action of~$\fg$
  and a~moment map in the sense of Lu $L =(L_+, L_-): V \to G^\ast$ such that $\Phi = L_+ L_-^{-1}$.
\end{Theorem}

\begin{Theorem}
    Let $V$ be a~finite-dimensional $G$-module which satisfies condition~\eqref{eq:key}, and which carries a~$G$-invariant
    nondegenerate symmetric bilinear form.
    Then, the corresponding Poisson space $(\Pi V, \pi)$ can be viewed as a Hamiltonian quasi-Poisson space
    with the same bivector, and with the group valued moment map $\Phi = \exp \circ \mu: V \to G$.
    Furthermore, the space $(\Pi V, \pi-r_{\Pi V})$ is Poisson, it carries an action of the Poisson-Lie group
    $G$ equipped with the standard Poisson bracket, and there is a~moment map in the sense
    of Lu $L =(L_+,  L_-): \Pi V \to G^\ast$ such that $\Phi = L_+ L_-^{-1}$.
\end{Theorem}

Our results may have an impact on the study of semiclassical counterparts of $q$-deformed Clifford algebras
and the corresponding quantum moment maps (see {\em e.g.} \cite{cubicDiracUqSL2}).
In more detail, consider a $G$-module $V$ satisfying condition~\eqref{eq:key}. The $G$-module $W=V \oplus V^*$ carries an invariant canonical scalar product $B$ which gives rise a constant Poisson bracket $\pi_B$ on $\Pi W$ and to the Clifford algebra ${\rm Cl}(W, B)$.
In order to define a Poisson-Lie deformation of $\pi_B$, we formulate the following problem:

\vskip 0.2cm

{\bf Problem:} Find a bivector $\pi_W$ on $\Pi W$ satisfying the following properties:
\begin{itemize}
\item
the action of $G$ on $(W,\pi_W)$ is Poisson;
\item 
$\pi_{\cW}$ descends to $r_{V}$ on~$V$, and to $r_{V^\ast}$ on~$V^\ast$;
\item
$\pi_W$ is compatible with the Poisson structure~$\pi_B$;
\item
$\pi_W$ is quadratic,
\end{itemize}

Assuming that such a bivector $\pi_W$ exists, the Poisson space 
$(\Pi W, \pi_B + \pi_W)$ would be a natural Poisson-Lie deformation of $(\Pi W, \pi_B)$, and its quantisation
would eventually provide a $q$-deformed Clifford algebra $\Cl_q(W)$.  Moreover, a~quantisation of
$(\Pi W, \pi_W)$ would provide a~$q$-deformed exterior algebra~$\twedge_qW$ which could be seen as a~fusion
of the $q$-deformed exterior algebras~$\twedge_qV$ and~$\twedge_qV^\ast$ defined
in~\cite{BerensteinZwicknagl2008}.
Such $q$-deformed exterior and Clifford algebras appear to be indispensable ingredients in the noncommutative differential geometry, see~\cite{HKdR,NicholsGrass,qSL2Diff}.
In this paper,  we obtain some partial results for~$V$ the defining representation of $\fg= \mathfrak{sl}_n$.

The structure of the paper is as follows: in Section 2, we collect the background material on Poisson and
quasi-Poisson actions, moment maps, dynamical $r$-matrices, and linear supermanifolds. In Section 3, we state
and prove our main results. In Section 4, we describe two simple examples related to the case of $\fg = \mathfrak{sl}_2$.
In Section~5, we discuss bivectors on the space $W=V \oplus V^*$.

\subsection*{Acknowledgements}
We would like to thank A.~\v{C}ap, D.~Leites, J.H.~Lu, P.~Pand\v{z}i\'c, and P.~\v{S}evera for stimulating discussion, and
the organizing committee of the conference ``Poisson 2024'' (July 8-12, 2024, Napoli, Italy) for an inspiring atmosphere during the meeting.
We are also grateful to the referee for constructive criticism and useful suggestions.
Research of A.~A. was supported in part by the grants 208235, 220040, and 236683, and by the National Center for Competence in Research (NCCR) SwissMAP of the Swiss National Science Foundation, and by the award of the Simons Foundation to the Hamilton Mathematics
Institute of the Trinity College Dublin under the program ``Targeted Grants to Institutes''.
Research of A.~K. was supported by the GA\v{C}R project 24-10887S.

\section{Preliminaries}
\subsection{Notation}

Let $G$ be a~complex simple Lie group and $\fg$ be its Lie algebra.
In what follows, we assume that we are working in the holomorphic category.
We denote by $\theta^L, \theta^R \in \Omega^1(G, \fg)$ the left- and right-invariant Maurer--\/Cartan forms on~$G$. For  $\beta\in\twedge\fg$ we denote by $\beta^L, \beta^R\in \Gamma(G, \twedge TG)$ the corresponding
left- and right-invariant multivector fields such that $\beta^L(e) = \beta^R(e) = \beta$.

The Lie algebra $\fg$ admits
a unique (up to a multiple) nondegenerate invariant symmetric bilinear form~$B_\fg$. Using~$B_\fg$, one can define the canonical
invariant Cartan 3-tensor on~$\fg$: 
\[
\phi(x,y,z) = \frac{1}{12} B_\fg(x, [y,z]).
\]
Using the bilinear form $B_\fg$ one can identify $\twedge \fg \cong \twedge \fg^*$, and by abuse of notation one can view $\phi$ as an element of $(\twedge^3 \fg)^\fg$.

For a $G$-manifold~$M$ and an element $x\in\fg$, the generating vector field of the induced  action of~$\fg$ is defined by
\[
  x_M(m) := \left.\frac{\dd}{\dd t}\right|_{t=0}\exp(-tx) \cdot m\qquad\text{for $m\in M$}.
\]
The Lie algebra homomorphism $\fg\to \fX(M)$ given by $x \mapsto x_M$ extends to a~map
\[
  \twedge\fg \to \Gamma(M, \twedge TM),
\]
preserving wedge products and Schouten brackets~$[\![-,-]\!]$. For $\alpha: M \to \twedge \fg$, we
denote by $\alpha_M$ the multivector field $\alpha_M(m) = (\alpha(m))_M$ for $m\in M$.
For $\pi \in \Gamma(M, \twedge^2 TM)$ a bivector on $M$, we denote by
$\pi^\sharp: T^*M \to TM$ the induced map between cotangent and tangent bundles
$\pi^\sharp_m(\alpha)(\beta) = \pi(\alpha,\beta)$ for $m \in M, \alpha, \beta \in T^*_mM$.

\subsection{Moment maps}

Recall that a~\emph{Hamiltonian Poisson $G$-manifold} is a triple $(M,\pi,\mu)$, consisting of a~$G$-manifold~$M$,
a $G$-invariant Poisson structure~$\pi$, and an~equivariant \emph{moment map} $\mu\colon M\to \fg^\ast$ satisfying the
condition
\[
  x_M = \pi^\sharp(\dd \langle \mu, x \rangle)  
  \qquad\text{for all $x\in\fg$}.
\]

In~\cite{AlekseevKosmannSchwarzbachMeinrenken2002}, a \emph{quasi-Poisson manifold} is defined as a~$G$-manifold~$M$ equipped with an invariant bivector
${\pi \in \Gamma(M, \twedge^2 TM)}$ such that
\[
  [\![\pi ,\pi]\!] = \phi_M.
\]
An $\Ad$-equivariant map $\Phi\colon M\to G$ is called a \emph{group valued moment map} if
\begin{equation}\label{eq:PhiCond}
  \pi^\sharp(\Phi^\ast B_\fg(x,\theta^L)) = \frac12((\id_\fg+\Ad_\Phi)x)_M\qquad\text{for all $x\in\fg$}.
\end{equation}
A~triple $(M, \pi, \Phi)$ is then called a \emph{Hamiltonian quasi-Poisson manifold}.

One says that the action of a~Poisson-Lie group $(G,\pi_G)$ on a~Poisson
$G$-manifold $(M,\pi)$ is a~\emph{Poisson action} if the action map
$G \times  M \to M$ is a Poisson map. 
If $G$ is connected, following~\cite{LuWeinstein1990,Lu1991} one can 
reformulate this condition as 
\[
  L_{x_M} \pi = -\delta(x)_{M},
\]
for all $x\in\fg$,  where $\delta\colon \fg \to \twedge^2\fg$ is the corresponding cobracket.

Let $G^\ast$ be the dual Poisson-Lie group of~$G$. A map $L\colon M \to G^\ast$  is called a~\emph{$G^\ast$-valued moment map} (or a moment map in the sense of Lu) if
\begin{equation}      \label{eq:momentLu}
   \pi^\sharp(L^\ast(\langle \theta^R, x \rangle)) = x_M
  \qquad\text{for all $x\in\fg$},
\end{equation}
where $\theta^R \in \Omega^1(G^\ast, \fg^*)$ is  the right-invariant Maurer-Cartan form on $G^\ast$.

Let $\fg=\fn_{+}\oplus\fh\oplus\fn_{-}$ be the
triangular decomposition of $\fg$. For $x \in \fg$, we denote by $x_+, x_0, x_-$ its components in $\fn_+, \fh, \fn_-$, respectively.
Let $H \subset G$ be a Cartan subgroup, $B_{\pm}=H N_\pm \subset G$ be the opposite Borel subgroups, and $\pr_{\pm}\colon B_{\pm} \to H$ be the natural projections.
For the Poisson-Lie group $G$ equipped with the Poisson structure defined by the standard $r$-matrix the dual Poisson-Lie group $G^\ast$ admits the following presentation:
\[
G^\ast = \left\{ (b_{+}, b_{-}) \in B_{+}\times B_{-} \mid \pr_{+}(b_{+}) \pr_{-}(b_{-})=e \right\} \subset G \times G,
\]
and the group $G$ can be viewed as the diagonal of $G \times G$. The corresponding Lie algebras are of the form
$$
\mathfrak{g} \cong \{(x,x); x \in \mathfrak{g}\}, \hskip 0.3cm
\mathfrak{g}^* \cong \{ (y,z); y \in \mathfrak{n}_+ \oplus \mathfrak{h},
z \in \mathfrak{n}_- \oplus \mathfrak{h}; y_0 + z_0=0\},
$$
and the pairing between $\mathfrak{g}$ and $\mathfrak{g}^*$ is realised using the invariant scalar product $B_\fg$:
$$
\langle (y,z), (x,x)\rangle = B_\fg(y-z, x).
$$

For the moment map $(L_+, L_-): M \to G^\ast$, the equation~\eqref{eq:momentLu} acquires the form
$$
\pi^\sharp(B_\fg(L_+^\ast \theta^R, x)) = \left(x_- + \frac{1}{2}x_0\right)_M, \hskip 0.3cm
\pi^\sharp(B_\fg(L_-^\ast  \theta^R, x )) = - \left(x_+ + \frac{1}{2}x_0\right)_M.
$$

We refer to~\cite{AlekseevMalkinMeinrenken1998} and~\cite{AlekseevKosmannSchwarzbach2000} for more details on
group valued moment maps.

\subsection{The dynamical $r$-matrix and group valued moment maps}\label{sec:Rdyn}
Let $e_a$ be an~orthonormal basis of~$\fg$ with respect to $B_\fg$. For a linear map $A\colon \fg\to\fg$, we define the corresponding matrix 
by $A(e_a) = \sum_b A_{ab}e_b$. Let $\varphi$ be a~meromorphic function of~$s\in\mathbb{C}$
given by 
\begin{equation}\label{eq:phi}
  \varphi(s) = \frac{1}{s} - \frac{1}{2}\coth\left(\frac{s}{2}\right)
  = - \frac{s}{12} + \frac{s^3}{720} - \frac{s^5}{30240} + O(s^7).
\end{equation}
Consider the map $T: \fg \to \twedge^2 \fg$  defined by the formula
\begin{equation}\label{eq:rDyn}
  T(x) = \frac{1}{2} \sum_{a,b} T_{ab}(x) e_a\wedge e_b,
\end{equation}
where $T_{ab}(x) = \varphi(\ad_x)_{ab}$.
By the results of~\cite{EtingofVarchenko1998}, $T$ is a~solution of the classical dynamical
Yang--\/Baxter equation
\[
  \mathrm{Cycl}_{abc}\left(\frac{\del T_{ab}}{\del e_c} + \sum_{kl}T_{ak}f_{akl}T_{lc}\right) = \frac14 \, f_{abc}
\]
where $\mathrm{Cycl}_{abc}$ denotes the sum over cyclic permutations and $f_{abc} = B_\fg(e_a,[e_b,e_c])$.

Let $(M,\pi)$ be a~Hamiltonian Poisson $G$-manifold with a moment map~$\mu\colon M\to\fg^\ast$, and
let $r_{\dyn}(\mu) = (\mu^\ast T)_M\in \Gamma(M, \twedge^2 TM)$,
where $\mu^*T: M \to \twedge^2 \fg$ is the pull-back of the function~$T$ under the moment map $\mu$.

By \cite[Theorem~7.1]{AlekseevKosmannSchwarzbachMeinrenken2002},
if the exponential map~$\exp\colon\fg\to G$ has the maximal rank on the image of~$\mu$,
then $\pi' = \pi + r_{\dyn}(\mu) $ defines a quasi-Poisson structure on~$M$ with the
group valued moment map $\Phi = \exp\circ \mu\colon M\to G$.

\subsection{Twist and $G^\ast$-valued moment maps}\label{sec:twist}
Let $\Delta^+(\fg)$ be the choice of positive roots of~$\fg$ given by the triangular decomposition of~$\fg$.
We assume that $B_\fg$ is such that the square of the length of the longest root is~2.
Let $\{E_\alpha\mid \alpha \in \Delta^+(\fg) \}$ be a~basis of root
vectors in~$\fn_+$ and $\{F_{\alpha}\mid \alpha \in \Delta^+(\fg) \}$ be the~basis of root vectors in~$\fn_{-}$ such that
$B_\fg(E_\alpha,F_\alpha)=1$.
Then the element
\[
  r = \frac12\sum_{\alpha\in\Delta^+(\fg)}  E_{\alpha}\wedge F_{\alpha}\in\twedge^2\fg
\]
is called the (antisymmetrised) standard $r$-matrix.
Note that
\[
   [\![ r , r ]\!] =  -\phi.
\]  

Let $H$ be a Cartan subgroup of~$G$, and $N_{+}$ and $N_{-}$ be nilpotent subgroups of~$G$
normalised by~$H$ such that $G_0 = N_{+}H N_{-}$ is an everywhere-dense subgroup of~$G$.
Then such a decomposition is known as the \emph{Gauss decomposition}; see~\cite[\S100]{ZhelobenkoBook}.
\ssbegin{Lemma}\label{lem:PhiFactorL}
  Let $(M,\pi, \Phi: M \to G)$ be a Hamiltonian quasi-Poisson $G$-manifold, 
and assume that ${\rm im}(\Phi)$ is contained in the domain of the Gauss decomposition $\Phi=L_+ L_-^{-1}$. Then, $\pi'=\pi-r_M$ is a Poisson bivector, and $L=(L_+, L_-): M \to G^*$ is the moment map in the sense of Lu.
\end{Lemma}

This seems to be a standard fact, but we have not found this formulation in the literature. Therefore, we include a sketch of the proof.

\begin{proof}
First, observe that
$$
[\![\pi' ,\pi']\!] = [\![\pi-r_M ,\pi-r_M]\!] = \phi_M - \phi_M =0.
$$
Hence, the bivector $\pi'$ is Poisson, as required.
Next, in equation \eqref{eq:PhiCond}, put $\Phi=L_+ L_-^{-1}$ and $x={\rm Ad}_{L_-} y$ to obtain
 $$
\pi^\sharp(B_\fg(L_+^* \theta^L - L_-^*\theta^L, y)) =
\frac{1}{2}({\rm Ad}_{L-}(y)+{\rm Ad}_{L+}(y))_M.
 $$
This implies
 $$
\pi^\sharp(B_\fg(L_+^* \theta^L, y)) = 
\begin{cases}
  \frac{1}{2}({\rm Ad}_{L-}(y)+{\rm Ad}_{L+}(y))_M & {\rm for} \,\,  y \in \mathfrak{n}_-, \\
  \frac{1}{4}({\rm Ad}_{L-}(y)+{\rm Ad}_{L+}(y))_M & {\rm for} \,\, y \in \mathfrak{h}.
\end{cases}
 $$
Here, the factor $\frac{1}{4}$ in the second line comes from the fact that
$(L_+^* \theta^L)_0 = -( L_-^*\theta^L)_0$.
 In order to compute the action of $r_M^\sharp$, we use the $G$-equivariance of $\Phi$ under the conjugation
 action of $G$ to obtain
$$
r_M^\sharp(L_+^* \theta^L) =
\frac{1}{2} \sum_{\alpha} \left(\iota(E_\alpha)(L_+^* \theta^L) (F_{\alpha})_M - \iota(F_{\alpha})(L_+^* \theta^L) (E_{\alpha})_M\right).
$$
For  components of the Gauss decomposition, the action by conjugation behaves as follows:
$$
\begin{array}{lll}
\iota(x_M) (L_+^* \theta^L) & = & -({\rm Ad}_{L_+^{-1}} x - {\rm Ad}_{L_-^{-1}} x)_+ - \frac{1}{2}({\rm Ad}_{L_+^{-1}} x - {\rm Ad}_{L_-^{-1}} x)_0, \\
\iota(x_M) (L_-^* \theta^L) & = & -({\rm Ad}_{L_-^{-1}} x - {\rm Ad}_{L_+^{-1}} x)_- - \frac{1}{2}({\rm Ad}_{L_-^{-1}} x - {\rm Ad}_{L_+^{-1}} x)_0.
\end{array}
$$
Here, the origin of the extra factors $\frac{1}{2}$ is the same as in the formula for $\pi^\sharp$ above.
Putting things together, we obtain
 $$
r_M^\sharp(B_\fg(L_+^* \theta^L, y)) = 
\begin{cases}
 \frac{1}{2}\left(
-({\rm Ad}_{L+}(y)-{\rm Ad}_{L-}(y))_- +({\rm Ad}_{L+}(y)-{\rm Ad}_{L-}(y))_+\right)_M &    y \in \mathfrak{n}_-, \\
\frac{1}{4}\left(
-({\rm Ad}_{L+}(y)-{\rm Ad}_{L-}(y))_- +({\rm Ad}_{L+}(y)-{\rm Ad}_{L-}(y))_+\right)_M & y \in \mathfrak{h}.
\end{cases}
 $$
Combining the expressions for $\pi^\sharp$ and $r_M^\sharp$, we get
 $$
\pi'^\sharp(B_\fg(L_+^* \theta^L, y)) = 
    \left(({\rm Ad}_{L_+} y)_- +\frac{1}{2} ({\rm Ad}_{L_+} y)_0 \right)_M. 
 $$
 By letting $z={\rm Ad}_{L_+} y$, we conclude 
 $$
\pi'^\sharp(L_+^*\langle \theta^R, z\rangle)= \pi'^\sharp(B_\fg(L_+^* \theta^R, z))=\pi'^\sharp(B_\fg(L_+^* \theta^L, y))=
\left(z_- + \frac{1}{2} z_0\right)_M.
 $$
The calculation for $L_-^* \theta^R$ is similar.
\end{proof}

\subsection{Exterior and symmetric algebras in the super case}
Let $V = V_\ev\oplus V_\od$ be a~supervector space. Recall that the exterior~$\mathbb{E}(V)$ and
symmetric~$\mathbb{S}(V)$ algebras of~$V$ are 
defined as
\[
  \mathbb{S}(V) = S(V_\ev) \otimes \twedge V_\od,\qquad
  \mathbb{E}(V) = \twedge V_\ev \otimes S(V_\od),
\]
where $\twedge W$ is the exterior algebra of a~vector space~$W$, and $S(W)$ is the symmetric
algebra of a~vector space~$W$.
Now, let $\cF$ be a~supercommutative algebra and $V$ be an $\cF$-module.
In the similar manner, one defines the exterior~$\mathbb{E}_{\cF}(V)$ and symmetric~$\mathbb{S}_{\cF}(V)$
algebras of~$V$. In what follows, we usually omit the subscript $\cF$.
We use $\wedge$ to denote the multiplication in the exterior algebra $\mathbb{E}(V)$ regardless of the
parity of elements, the multiplication in the symmetric algebra~$\mathbb{S}(V)$ is denoted by juxtaposition.

\subsection{Linear supermanifolds}
For an introduction to supergeometry, we refer to~\cite{Berezin,QFSsuper,Leites80,SoS1rus},
for the basics of complex supergeometry, see~\cite{ManinGaugeFields,KostantGraded} and review in~\cite{Leites2023}.
Complex supermanifolds
are ringed spaces, i.e., pairs $\cM := (M, \cO_{\cM})$, where $M$ is an $m$-dimensional complex manifold and
$\cO_{\cM}$ is the structure sheaf of~$\cM$, locally isomorphic to $\cO_U\otimes\twedge(n)$, where $\cO_U$ is
the sheaf of holomorphic function on a~domain~$U$ in~$M$, and $\twedge(n)$ is the Grassmann algebra with $n$
generators.
A morphism of supermanifolds $f\colon \cM\to \cN$ is a~morphism of the corresponding ringed spaces (for
details, see~\cite[\S2.3.1]{ManinSchemes}) that preserves the parity of the sections.
A linear supermanifold corresponding to a~super vector space $V=V_\ev\oplus V_\od$ is defined by
\[
  \cV = (V_\ev, \cO_{V_\ev}\otimes \twedge V_\od^\ast),
\]
where $\cO_{V_\ev}$ is the sheaf of holomorphic functions on~$V_\ev$.

Since we work with linear supermanifolds, without loss of generality, we can replace the structure
sheaf of~$\cV$ by the algebra of holomorphic function on~$\cV$ which we denote by~$\cF$.
We consider the subalgebra $\cF_{\poly}$ of polynomial functions on a~linear supermanifold~$\cV$ which can be
identified with~$\mathbb{S}(V^\ast)$. The algebra $\cF_{\poly}$ has a~natural $\Zee$-grading given by
\begin{equation}\label{eq:FpolyGrad}
  \cF_{\poly} = \bigoplus_{i} \cF_{\poly}^i,\qquad
  \cF_{\poly}^i := \mathbb{S}^i(V^\ast)
  \cong \bigoplus_{a+b=i} S^a(V_\ev^\ast)\otimes \twedge^bV_\od^\ast.
\end{equation}

The vector fields on~$\cV$ are elements of the free $\cF$-module $\fX(\cV) := \Der\cF$.
Multivector fields are elements of the corresponding exterior algebra~$\mathbb{E}_{\cF}\fX(\cV)$. In
particular, bivectors are elements of~$\mathbb{E}_{\cF}^2\fX(\cV)$.

In what follows, we assume that~$V$ is a~finite-dimensional super vector space equipped with a~structure of
a~$G$-module and a~non-degenerate invariant skew-symmetric even bilinear form~$B$. Often we identify~$V$ and
$V^\ast$ using~$B$. We consider the following two cases:

1) $\dim V_\od = 0$, and $B$ defines a skew-symmetric bilinear form on $V_\ev$. We call such modules
even, and we have  \(\cF_{\poly} = \mathbb{S}(V^\ast) \cong S(V_\ev^\ast)\). 
Using the isomorphism $V^\ast \cong V$ induced by $B$, we can view elements of a basis
$\{ v_i \}$ of $V$ as generators of~$\cF_{\poly}$.
The corresponding basis~$\del_i$ of $\fX(\cV)$ is defined by $\del_i(v_j) = \delta_{i,j}$.

2) $\dim V_\ev = 0$, and $B$ is a symmetric bilinear form on~$V_\od$. We call such modules odd, and we have \(\cF=\cF_{\poly} = \mathbb{S}(V^\ast) \cong \twedge V_\od^\ast\) . As before, $B$ induces an isomorphism $V^\ast \cong V$, and we can view elements of a basis $\{ \xi_i\}$ of $V$ as generators of~$\cF_{\poly}$.
The corresponding basis~$\del_i$ of $\fX(\cV)$ is defined by $\del_i(\xi_j) = \delta_{i,j}$.

In what follows, we will consider only multivector fields with polynomial coefficients.
Therefore, we have the following isomorphism of $G$-modules
\[
  \fX_{\poly}(\cV) \cong  \cF_{\poly}\otimes V.
\]
The $\Zee$-grading~\eqref{eq:FpolyGrad} induces a $\Zee$-grading on the Lie
superalgebra~$\fX_{\poly}(\cV) = \bigoplus_{i={-1}}^{\infty}\fX^i_{\poly}(\cV)$ of vector fields
with polynomial coefficients by setting $\deg\del_i = -1$.
This grading is extended to a $\Zee^2$-grading of the multivector fields with polynomial coefficients
\begin{equation}\label{eq:MVectGrad}
  \mathbb{E}_{\cF}\fX_{\poly}(\cV) = \bigoplus\nolimits_{i,j} \mathbb{E}^{i,j} \fX_{\poly}(\cV),
\end{equation}
where $i$ denotes the degree with respect to the standard grading of the exterior algebra and $j$ denotes the
degree induced by the standard grading of~$\fX_{\poly}(\cV)$.

\subsection{Moment maps in the super case}
The notions from (quasi-)Poisson geometry described above have straightforward generalisations to the
super setting. We  associate to~$\fg^\ast$, $G$ and~$G^\ast$ the corresponding even supermanifold. Then, the moment maps
$\mu\colon \cM\to \fg^\ast$, $\Phi\colon \cM\to G$, $L\colon \cM\to G^\ast$ are defined as morphisms of
supermanifolds satisfying the same conditions as in the classical case.
Another interpretation of the moment maps is to consider them in terms of $\cM$-points of $\fg^\ast$, $G$, and
$G^\ast$, see~\cite[\S\S2.8--2.9]{QFSsuper}.

\section{Moment maps and classical $r$-matrices}

\subsection{Poisson brackets and bilinear invariants}
Let $V$ be a~finite dimensional $\fg$-module which admits a~non-degenerate $\fg$-invariant symmetric or symplectic bilinear form~$B$. Then, the formula
\begin{equation}\label{eq:brB}
  \{v,w\}_B = 2B(v,w)\quad v,w\in V
\end{equation}
defines a~$\fg$-equivariant Poisson bracket on~$\twedge V$ if $B$ is symmetric, and on~$S(V)$ if $B$ is skew-symmetric.
Therefore, we obtain a Poisson structure on the linear supermanifold~$\cV$ associated with the super vector space~$V$ if $V$
is even, and with the super vector space $\Pi V$ if $V$ is odd. We denote by $\pi_B\in(\mathbb{E}^2\fX(\cV))^\fg$ the corresponding bivector on~$\cV$.
As it was shown in~\cite{KostantSternberg}, the action of~$\fg$ on~$\cV$ is Hamiltonian, and the corresponding moment map is given by
\begin{equation}\label{eq:mu}
  \mu  = \frac14\sum_{i=1}^{\dim\fg} \sum_{j=1}^{\dim V} (e_i.v_j)v_j^\ast\otimes f_i \in
  \cF^2_{\poly}\otimes \fg^\ast,
\end{equation}
where $e_i$ is a basis of~$\fg$ and $f_i$ is the corresponding dual basis of~$\fg^\ast$,
$v_j$ is a~basis of~$V$ and $v_j^\ast$ is the corresponding dual basis with respect to~$B$.

Nondegenerate invariant bilinear forms on finite-dimensional irreducible $\fg$-modules were classified
in~\cite{Malcev1944sub,Dynkin1950some,Dynkin1951auto}, see also a review in~\cite[Theorem~0.20]{Dynkin1952}.

\subsection{The $r$-matrix Poisson bracket}
Let $r\in\twedge^2\fg$ be the (antisymmetrised) classical $r$-matrix of~$\fg$.
For a $\fg$-module~$V$, 
we define a~quadratic bracket $\{-,-\}_r$ on the symmetric and exterior algebras of~$V$:
\begin{equation}\label{eq:brR}
  \{a,b\}_r = ({\rm mult} \circ r_{V}) (a \otimes b).
\end{equation}
Here ${\rm mult}$ is the product in the symmetric or exterior algebra,
and $r_{V} (a \otimes b)$ is the action of a bivector on a pair of functions.
The bracket \eqref{eq:brR} is by construction skew-symmetric and satisfies the Leibniz rule. By \cite[Main
Theorem~1.1]{Zwicknagl2009}, it is Poisson if and only if $\Hom_\fg(\twedge^3V,S^3(V)) =0$,
see also~\cite[Proposition~5.2]{GoodearlYakimov2009} for a~partial list of such $\fg$-modules, and \cite{BerensteinZwicknagl2008} for related quantisation problems.

Denote by  $\varpi_i$  the $i$th fundamental weight of~$\fg$ (the labeling of the fundamental weights follows  Table~1 on p.~293 in~\cite{VinbergOnishchik}).

\ssbegin{Lemma}\label{lem:LIST}
Let $(\fg,V)$ be a~pair of a Lie algebra~$\fg$ and a~selfdual irreducible finite-dimensional $\fg$-module such that $\{-,-\}_r$ is a~Poisson bracket on~$\cV$. Then,

1) if~$B$ is symmetric, then the pair $(\fg,V)$ belongs to the following list:
$(\fsl_2,V_{2\varpi_1})$, $(\fsl_4,V_{\varpi_2})\cong(\fo_6,V_{\varpi_1})$, $(\fo_8,V_{\varpi_3})$,  $(\fo_8,V_{\varpi_4})$,  $(\fo_n,V_{\varpi_1})$

2) if~$B$ is skew-symmetric, then the pair $(\fg,V)$ belongs to the following list:
$(\fsl_2,V_{\varpi_1})$, $(\fsp_{2n},V_{\varpi_1})$,
$(\fo_5,V_{\varpi_2}) \cong (\fsp_4,V_{\varpi_1})$.
\end{Lemma}

\begin{proof}
  The list of pairs $(\fg,V)$ such that $\{-,-\}_r$ is a~Poisson bracket on~$\cV$ is given in~\cite[Main
  Theorem~3.12]{Zwicknagl2009}:
  \begin{enumerate}
  \item $(\fsl_n,V_\lambda)$, where  $\lambda \in
    \{\varpi_1,2\varpi_1,\varpi_2,\varpi_{n-2},2\varpi_{n-1},\varpi_{n-1} \}$,
  \item $(\fo_n,V_{\varpi_1})$, $(\fo_5,V_{\varpi_2})$, $(\fo_8,V_{\varpi_3})$, $(\fo_8,V_{\varpi_4})$,
    $(\fo_{10},V_{\varpi_4})$, $(\fo_{10},V_{\varpi_5})$,
  \item $(\fsp_{2n},V_{\varpi_1})$, $(\fsp_4, V_{\varpi_2})$,
  \item $(\fe_6,V_{\varpi_1})$, $(\fe_6,V_{\varpi_5})$.
  \end{enumerate}

  As it was shown in~\cite{Malcev1944sub,Dynkin1950some,Dynkin1951auto}, an irreducible finite-dimensional $\fg$-module~$V$
  admits a~nondegenerate invariant bilinear form if

  \begin{itemize}
  \item for $\fsl_n$, $\fo_{4k+2}$, $\fe_6$, the marks of the highest weight of~$V$ on the Dynkin diagram are
    distributed symmetrically;
  \item for $\fsp_{2n}$, $\fo_{2n+1}$, $\fo_{4n}$, such bilinear form always exists.
  \end{itemize}

  Now the claim follows from comparing the lists.
\end{proof}

The following theorem shows that odd and even modules listed in Lemma~\ref{lem:LIST} admit the following
remarkable properties.

\ssbegin{Theorem}      \label{thm:main}
Let $V$ be a~finite-dimensional irreducible $\fg$-module such that
\begin{enumerate}\renewcommand{\theenumi}{\alph{enumi}}
\item $V$ admits a~nondegenerate invariant (symmetric or symplectic) bilinear form $B$,
\item $\Hom_\fg(\twedge^3 V, S^3(V)) =0$.
\end{enumerate}
Then,
\begin{enumerate}
\item 
the space $(V, \pi_B)$ (or $(\Pi V, \pi_B)$) is a Hamiltonian quasi-Poisson space with group valued moment map $\Phi=\exp \circ \mu$;
\item $r_{\dyn}(\mu) = 0$;
\item 
  the space $(\cV, \pi_B - r_{\cV})$ is Poisson and carries an action of the Poisson-Lie group~$G$ equipped with the standard
  Poisson bracket;
\item 
for $B$ symmetric, the $G$-action on $(\Pi  V, \pi - r_{\cV}))$ admits a~moment map in the sense of Lu $L = (L_{+}, L_{-}): \cV \to G^\ast$, where $\Phi = L_{+}L_{-}^{-1}$;
  \item 
  for $B$ symplectic, 
  the pre-image under $\Phi$ of the domain of the Gauss decomposition carries a Hamiltonian $\fg$-actions and a moment
  map in the sense of Lu $L = (L_{+}, L_{-}):  \cV \to G^\ast$, where $\Phi = L_{+}L_{-}^{-1}$.
\end{enumerate}
\end{Theorem}

The proof is divided into several parts.

\ssbegin{Lemma}\label{lem:quasi}
Let $V$ be an even or odd $\fg$-module such that
$\Hom_\fg(S^3(V),\twedge^3V)=0$, and let $\cV$ be the corresponding linear supermanifold. Then, $\phi_{\cV} = 0$
and every $G$-invariant Poisson structure on~$\cV$ can be viewed as a quasi-Poisson structure.
\end{Lemma}
\begin{proof}
For $x\in\fg$ the Lie algebra homomorphism $\fg\to\fX(\cV)$ is given by
\begin{align*}
  x\mapsto x_{\cV} =&{} \sum_i (x.\xi_i)\del_i \in \fX_{\poly}^0(\cV) \cong \twedge^1V^\ast\otimes S^1(V)
  &{}& \text{for $V$ odd},\\
  x\mapsto x_{\cV} =&{} \sum_i (x.v_i)\del_i \in \fX_{\poly}^0(\cV) \cong S^1(V^\ast)\otimes \twedge^1V
  &{}& \text{for $V$  even},
\end{align*}
Therefore, the image of the Cartan element satisfies the following condition
\begin{align*}
  \phi_\cV \in {}& \mathbb{E}^{3,0}\fX_{\poly}(\cV)^\fg \cong \left(\twedge^3 V^\ast \otimes S^3(V)\right)^\fg & & \text{for $V$  odd},\\
  \phi_\cV \in {}& \mathbb{E}^{3,0}\fX_{\poly}(\cV)^\fg \cong \left(S^3(V^\ast) \otimes \twedge^3 V\right)^\fg & & \text{for $V$ even}.
\end{align*}
Using $B$, we identify $\twedge V \cong \twedge V^\ast$ and $S(V) \cong S(V^\ast)$.
In  both cases, we have the following isomorphism of $\fg$-modules
\[
  \phi_{\cV} \in \mathbb{E}^{3,0}\fX_{\poly}(\cV)^\fg \simeq \Hom_\fg(\twedge^3V, S^3(V)).
\]
In particular, if $\Hom_\fg(\twedge^3 V, S^3(V)) =0$, then $\phi_{\cV} = 0$. This proves the claim.
\end{proof}

\ssbegin{Lemma}\label{lem:Rdyn}
Let $V$ be a $\fg$-module listed in Lemma~\ref{lem:LIST}, then $r_{\dyn}(\mu)=0$.
\end{Lemma}
\begin{proof}
  First, we will show that for all $(\fg,V)$ listed in Lemma~\ref{lem:LIST} we have
  \begin{equation}\label{eq:admuk}
    \ad_\mu^3 = 0.
  \end{equation}
  Using the isomorphism $\fg\cong\fg^\ast$, and the facts that the moment map $\mu$ is quadratic and $G$-equivariant, we obtain
  \[
    \ad_\mu^k \in \Hom_\fg(\fg,\cF^{2k}_{\poly}\otimes\fg) \cong \Hom_\fg(\fg\otimes\fg, \cF^{2k}_{\poly}).
  \]
  Therefore, the equality~\eqref{eq:admuk} would follow from $\Hom_\fg(\fg\otimes\fg, \cF^{2k}_{\poly}) = 0$
  for $k=3$. In order to prove this fact, we use Table~5 on p.~299 in~\cite{VinbergOnishchik} or direct computations.

  $\bullet$ For $(\fg,V)=(\fsl_2,V_{\varpi_1})$, we have
  \[
    \cF^{6}_{\poly} \cong S^6(V) \cong V_{6\varpi_1},\qquad
    \fg\otimes\fg \cong V_{4\varpi_1}\oplus V_{2\varpi_1}\oplus \Cee,
  \]
  which implies $\Hom_\fg(\fg\otimes\fg, \cF^{6}_{\poly}) = 0$ and $\ad_\mu^3 = 0$.

  $\bullet$ For $(\fg,V)=(\fsp_{2n},V_{\varpi_1})$, we have
  \[
    \cF^{6}_{\poly} \cong S^6(V) \cong V_{6\varpi_1},\qquad
    \fg\otimes\fg \cong V_{2\varpi_1+\varpi_2}\oplus V_{2\varpi_1}
    \oplus V_{4\varpi_1}\oplus V_{2\varpi_2}\oplus V_{\varpi_2}\oplus \Cee,
  \]
  which implies $\Hom_\fg(\fg\otimes\fg, \cF^{6}_{\poly}) = 0$ and $\ad_\mu^3=0$.

  $\bullet$ For $(\fg,V)=(\fsl_{2},V_{2\varpi_1})$, we have ${\rm dim} V_{2\varpi_1} =3$ which implies
  $\cF^{4}_{\poly}=\twedge^4V_{2\varpi_1}=0, \Hom_\fg(\fg\otimes\fg, \cF^{4}_{\poly}) = 0$, and $\ad_\mu^2=0$.

  $\bullet$ For $(\fg,V)=(\fo_n,V_{\varpi_1})$,
  let $\xi_1,\ldots,\xi_n$ be an orthonormal basis in~$V=V_{\varpi_1}$ and $E_{i,j}  - E_{j,i}$, $i<j$ be a basis  in~$\fo_n$.
We have that
\[
  \frac14 \sum_a ((E_{i,j}-E_{j,i}).\xi_a)\xi_a^\ast
  = \frac14\left(\xi_j\xi_i - \xi_j\xi_i\right) = \frac12 \xi_i\xi_j.
\]
Let $B_\fg(x,y) = \tr_V(x y)$ for $x,y\in\fo_n$, then
\[
  (E_{i,j} - E_{j,i})^\ast = - \frac12(E_{i,j}-E_{j,i}).
\]
Set $\mu_{i,j} = -\frac14 \xi_i\xi_j$, $e_{i,j} =  E_{i,j} - E_{j,i}$, then we get that
\[
  \mu = -\frac14 \sum_{i<j} \xi_i\xi_j\otimes (E_{i,j} - E_{j,i}) = \sum _{i<j} \mu_{i,j}\otimes e_{i,j}.
\]
We have that
\begin{align*}
  \ad_{\mu}^3(e_{a,b})  ={}
  & \sum_{x<y;i<j;k<l} \mu_{i,j}\mu_{x,y}\mu_{k,l}\otimes [e_{i,j}, [e_{x,y}, [e_{k,l}, e_{a,b}]]].
\end{align*}
Since $\xi_i$ are odd, $\mu_{x,y}\mu_{i,j}\mu_{k,l}$ is nonzero only when all indexes $i,j,x,y,k,l$ are distinct.

It is easy to see that for any $a,b$ and distinct $i,j,x,y,k,l$ we have that
\[
  [e_{i,j}, [e_{x,y}, [e_{k,l}, e_{a,b}]]] = 0.
\]
Therefore, we get that $\ad_\mu^3=0$.

Due to the symmetry of the Dynkin diagram of~$\fo_8$, similar computations show that $\ad_\mu^3=0$ for $(\fo_8,V_{\varpi_3})$,
$(\fo_8, V_{\varpi_4})$, and $(\fo_4,V_{\varpi_1}) = (\fsl_4,V_{\varpi_2})$.

\medskip
To summarise, for all modules from the list in Lemma~\ref{lem:LIST}
  we have
  \[
    \varphi(\ad_\mu) = -\frac{1}{12}\ad_\mu,
  \]
  
  Next, we will show that the corresponding terms in $r_{\dyn}(\mu)$ vanish.
  Using~\eqref{eq:phi} and~\eqref{eq:rDyn}, and the fact that $\ad_\mu^3=0$
  in all cases of interest, we obtain
  \begin{equation}\label{eq:rDynSer}
    r_{\dyn}(\mu) = \left(-\frac{1}{12}\sum_{ab} (\ad_\mu)_{ab} e_a\wedge e_b
    \right)_{\cV} \in 
    (\mathbb{E}^{2,4}\fX_{\poly}(\cV))^\fg.
  \end{equation}
  Note that
  \begin{align*}
    (\mathbb{E}^{2,4}\fX_{\poly}(\cV))^\fg \cong {}
    & \Hom_\fg(S^{4}(V), \twedge^2 V) & &\text{for $V$ even}, \\
    (\mathbb{E}^{2,4}\fX_{\poly}(\cV))^\fg \cong {}
    & \Hom_\fg(\twedge^{4}V, S^2(V)) & &\text{for $V$  odd}.
  \end{align*}
  In particular, we have
  \begin{equation}\label{eq:muSquareDyn}
    \sum_{a,b} ((\ad^{2k-1}_\mu)_{ab} e_a\wedge e_b)_{\cV}\in (\mathbb{E}^{2,4k}\fX_{\poly}(\cV))^\fg
    \cong \begin{cases}
      \Hom_\fg(S^{4k}(V), \twedge^2 V) & \text{for $V$  even}, \\
      \Hom_\fg(\twedge^{4k}V, S^2(V))  &\text{for $V$  odd}.
    \end{cases}
  \end{equation}
  Using Table~5 on p.~299 in~\cite{VinbergOnishchik}, for each pair $(\fg,V)$ listed in Lemma~\ref{lem:LIST} we
  show that the expression \eqref{eq:muSquareDyn} vanishes for $k=1$  which implies $r_{\dyn}(\mu) = 0$.

  We start with even modules.
  
  $\bullet$ For $(\fg,V) = (\fsl_2,V_{\varpi_1})$, we get the following isomorphisms of $\fg$-modules
  \[
    S^4(V) \cong V_{4\varpi}\quad\text{and}\quad
    \twedge^2 V \cong \mathbb{C}
  \]
  which implies $\Hom_\fg(S^4(V), \twedge^2 V) = 0$.

  $\bullet$ For $(\fg,V) = (\fsp_{2n},V_{\varpi_1})$, we get the following isomorphisms of $\fg$-modules
  \[
    S^4(V) \cong V_{4\varpi_1}\quad\text{and}\quad
    \twedge^2 V \cong V_{\varpi_2}\oplus\mathbb{C}
  \]
  which implies $\Hom_\fg(S^4(V), \twedge^2 V) = 0$.

  $\bullet$ For $(\fg,V) = (\fsp_4,V_{\varpi_1}) \cong (\fo_5,V_{\varpi_2})$, we get the following isomorphisms of $\fg$-modules
  \[
    S^4(V) \cong V_{4\varpi},\quad\text{and}\quad
    \twedge^2 V \cong \mathbb{C}
  \]
  which implies $\Hom_\fg(S^4(V), \twedge^2 V) = 0$.

  Next, we consider odd modules listed in Lemma~\ref{lem:LIST}.
  
  $\bullet$ For $(\fg,V) = (\fsl_2,V_{2\varpi_1})$, we have $\twedge^4V = 0$, and $\Hom_\fg(\twedge^4V, S^2( V)) = 0$.

  $\bullet$ For $(\fg,V) = (\fsl_4,V_{\varpi_2})$, we get the following isomorphisms of $\fg$-modules
  \[
    \twedge^4V \cong V_{\varpi_1+\varpi_3},\qquad
    S^2(V) \cong V_{2\varpi_2} \oplus \mathbb{C},
  \]
  which implies $\Hom_\fg(\twedge^4V, S^2( V)) = 0$. 

  $\bullet$ For $(\fg,V) = (\fo_8,V_{\varpi_3})$, we get the following isomorphisms of $\fg$-modules
  \[
    \twedge^4V \cong  V_{2\varpi_1}\oplus V_{2\varpi_4},\qquad
    S^2(V) \cong V_{2\varpi_3}\oplus \mathbb{C},
  \]
  which implies $\Hom_\fg(\twedge^4V, S^2( V)) = 0$.
  
  $\bullet$ For $(\fg,V) = (\fo_8,V_{\varpi_4})$, we get the following isomorphisms of $\fg$-modules
  \[
    \twedge^4V \cong  V_{2\varpi_1}\oplus V_{2\varpi_3},\qquad
    S^2(V) \cong V_{2\varpi_4}\oplus \Cee,
  \]
  which implies that $\Hom_\fg(\twedge^4V, S^2( V)) = 0$. 
  
  $\bullet$ For $(\fg,V) = (\fo_n,V_{\varpi_1})$, we get the following isomorphisms of $\fg$-modules
  \begin{subequations}\label{eq:OnE4V1}
    \begin{align}
      \twedge^4V \cong {}&V_{\varpi_1} & & \text{for $\fo_5$},\\
      \twedge^4V \cong {}& V_{\varpi_2+\varpi_3} & & \text{for $\fo_6$},
    & \twedge^4V \cong {}& V_{2\varpi_3} & & \text{for $\fo_7$},\\
      \twedge^4V \cong {}& V_{2\varpi_3}\oplus V_{2\varpi_4} & & \text{for $\fo_{8}$},
        & \twedge^4V \cong {}& V_{2\pi_4} & & \text{for $\fo_9$},\\
      \twedge^4V \cong {}& V_{\varpi_4+\varpi_5} & & \text{for $\fo_{10}$},
    & \twedge^4V \cong {}& V_{\varpi_4} & & \text{for $\fo_n$, $n>10$},
    \end{align}
  \end{subequations}
  and 
  \[
    S^2(V) \cong V_{2\pi_1}\oplus \mathbb{C},
  \]
  which implies  $\Hom_\fg(\twedge^4V, S^2( V)) = 0$.

  This implies that $r_{\dyn}(\mu) = 0$ for all pairs $(\fg,V)$ listed in  Lemma~\ref{lem:LIST}.
\end{proof}

\ssbegin{Lemma}\label{lem:PhiMoment}
Let $(\fg,V)$ be a~pair listed in Lemma~\ref{lem:LIST}, $\mu$ be the corresponding moment map,
and  $W=V_{\varpi_1}$ be the standard representation of~$\fg$. By identifying
$\mu \in \cF\otimes \End(W)$ and $\exp\circ\mu \in \cF\otimes \End(W)$, we get that

1) $\exp\circ\mu = 1 +\mu$ for $(\fsl_2,V_{\varpi_1})$, $(\fsp_{2n},V_{\varpi_1})$, $(\fsl_2,V_{2\varpi_1})$,
$(\fo_5,V_{\varpi_1}$, $(\fo_n,V_{\varpi_1})$, where $n>6$;

2) $\exp\circ\mu = 1 + \mu + \frac12\mu^2$ for $(\fo_6,V_{\varpi_1})\cong(\fsl_4,V_{\varpi_2})$;

3) $\exp\circ\mu = 1 + \mu + \tfrac12\mu^2 + \tfrac16\mu^3 + \tfrac1{24}\mu^4$ for $(\fo_8,V_{\varpi_3})$ and $(\fo_8,V_{\varpi_4})$.

\noindent Moreover, $\Phi$ has maximal rank.
\end{Lemma}
\begin{proof}
  To compute $\exp\circ\mu$ we consider $\mu$ as an element of $(\cF_{\poly}^2\otimes\End(V))^\fg$.
  Let $W=V_{\varpi_1}$ be the standard representation for~$\fg$.
  We have 
  \[
    \exp\circ\mu \in\bigoplus_{k=0}^{\infty} (\cF_{\poly}^{2k}\otimes\End(W))^\fg
    \cong \left(\cF_{\poly}^{2k}\otimes (W\otimes W^\ast)\right)^\fg.
  \]
  Using Table~5 on p.~299 in~\cite{VinbergOnishchik}, for each pair $(\fg,V)$ listed in Lemma~\ref{lem:LIST} we
  show that the space $(\cF_{\poly}^{2k}\otimes\End(W))^\fg$ vanishes for certain $k$.

  
  $\bullet$ For $(\fg,V)=(\fsl_2,V_{\varpi_1})$ we have that $W=V_{\varpi_1}$ and
  \[
    \left(\cF_{\poly}^{2k}\otimes (W\otimes W^\ast)\right)^\fg
    \cong  \left(S^{2k}(V) \otimes (V_{2\pi}\oplus \mathbb{C})\right)^\fg
    \cong  \left(V_{2k\pi} \otimes (V_{2\pi}\oplus \mathbb{C})\right)^\fg    
  \]
  which is zero if $k>1$. Therefore, $\exp\circ\mu = \id + \mu$.
  
  $\bullet$ For $(\fg,V)=(\fsp_{2n},V_{\varpi_1})$ we have that $W=V_{\varpi_1}$ and
  \[
    \left(\cF_{\poly}^{2k}\otimes (W\otimes W^\ast)\right)^\fg
    \cong  \left(S^{2k}(V) \otimes (V_{2\varpi_1}\oplus V_{\varpi_2} \oplus \mathbb{C})\right)^\fg
    \cong \left(V_{2k\varpi_1}\otimes (V_{2\varpi_1}\oplus V_{\varpi_2} \oplus \mathbb{C})  \right)^\fg
  \]
  which is zero if $k>1$. Therefore, $\exp\circ\mu = \id + \mu$.

  $\bullet$ For $(\fg,V) = (\fsp_4,V_{\varpi_1}) \cong (\fo_5,V_{\varpi_2})$,
  the calculation is the same as for $(\fsp_{2n},V_{\varpi_1})$.
  Therefore, $\exp\circ\mu = \id + \mu$.

  $\bullet$ For $(\fg,V)=(\fsl_2,V_{2\varpi_1})$, $W=V_{\varpi_1}$ we have that
  \[
    \left(\cF_{\poly}^{2k}\otimes (W\otimes W^\ast)\right)^\fg =  0\qquad\text{for $k>1$},
  \]
  since $\cF_{\poly}^{2k}\cong\twedge^{2k}V_{2\varpi} = 0$ for $k>1$. Therefore, $\exp\circ\mu = \id + \mu$.

  $\bullet$ For $(\fg,V) = (\fo_8,V_{\varpi_3})$, $W=V_{\varpi_1}$ we have that
  \[
    \left(\cF_{\poly}^{2k}\otimes (W\otimes W^\ast)\right)^\fg
    \cong  \left(\twedge^{2k}V_{\varpi_3} \otimes (V_{\varpi_2}\oplus V_{2\varpi_1}\oplus\mathbb{C})\right)^\fg    
  \]
  We have that
  \[
    \twedge^2V_{\varpi_3}\cong V_{\varpi_2},\quad
    \twedge^4V_{\varpi_3}\cong V_{2\varpi_1}\oplus V_{2\varpi_4},\quad
    \twedge^6V_{\varpi_3} \cong V_{\varpi_2},\quad
    \twedge^8V_{\varpi_3}  \cong \Cee.
  \]
  Since $\dim V_{\varpi_3} = 8$, we get that $\exp\circ\mu = 1 + \mu + \tfrac12\mu^2 + \tfrac16\mu^3 + \tfrac1{24}\mu^4$.

  $\bullet$ For $(\fg,V) = (\fo_8,V_{\varpi_4})$, $W=V_{\varpi_1}$
  \[
    \left(\cF_{\poly}^{2k}\otimes (W\otimes W^\ast)\right)^\fg
    \cong  \left(\twedge^{2k}V_{\varpi_4} \otimes (V_{\varpi_2}\oplus V_{2\varpi_1}\oplus\mathbb{C})\right)^\fg    
  \]
  We have that
  \[
    \twedge^2V_{\varpi_4}\cong V_{\varpi_2},\quad
    \twedge^4V_{\varpi_4}\cong V_{2\varpi_1}\oplus V_{2\varpi_3},\quad
    \twedge^6V_{\varpi_4} \cong V_{\varpi_2},\quad
    \twedge^8V_{\varpi_4}  \cong \Cee.
  \]
  Since $\dim V_{\varpi_4} = 8$, we get that $\exp\circ\mu = 1 + \mu + \tfrac12\mu^2 + \tfrac16\mu^3 + \tfrac1{24}\mu^4$.

  $\bullet$ For $(\fg,V) = (\fo_5,V_{\varpi_1})$,  $W=V_{\varpi_1}$ we have that
  \[
    \left(\cF_{\poly}^{2k}\otimes (W\otimes W^\ast)\right)^\fg
    \cong  \left(\twedge^{2k}V_{\varpi_1} \otimes (V_{2\varpi_1}\oplus V_{2\varpi_2}\oplus \Cee)\right)^\fg    
  \]  
  Since $\dim V_{\varpi_1} = 5$, $\twedge^{2k}V_{\varpi_1} = 0$ for $k>2$. We have that
  \[
    \twedge^2V_{\varpi_1} \cong V_{2\varpi_2},\qquad \twedge^4 V_{\varpi_1} \cong V_{\varpi_1}.
  \]
  Therefore, $\exp\circ\mu = 1 + \mu$.
  
  $\bullet$ For $(\fg,V) = (\fo_6,V_{\varpi_1})$,  $W=V_{\varpi_1}$ we have that
  \begin{equation}\label{eq:isoPhiO6}
    \left(\cF_{\poly}^{2k}\otimes (W\otimes W^\ast)\right)^\fg
    \cong  \left(\twedge^{2k}V_{\varpi_1} \otimes (V_{2\varpi_1}\oplus V_{\varpi_2+\varpi_3}\oplus \Cee)\right)^\fg    
  \end{equation}
  We have that
  \[
    \twedge^2V_{\varpi_1}\cong \twedge^4V_{\varpi_1}\cong V_{\varpi_2+\varpi_3},\qquad
    \twedge^6V_{\varpi_1}\cong \mathbb{C}.
  \]
  Since $\dim V_{\varpi_1} = 6$, \eqref{eq:isoPhiO6} vanishes for $k>3$. 
  Therefore, $\exp\circ\mu = 1 + \mu + \frac12\mu^2 + \frac16\mu^3$.
  
  Let us note that if we consider the pair~$(\fo_6,V_{\varpi_1})$ as $(\fsl_4,V_{\varpi_2})$ (hence $W\cong
  V_{\varpi_2}$ as $\fo_6$-modules), then we get that
  \[
    \left(\cF_{\poly}^{2k}\otimes (W\otimes W^\ast)\right)^\fg
    \cong  \left(\twedge^{2k}V_{\varpi_2} \otimes (V_{\varpi_1+\varpi_3}\oplus \mathbb{C})\right)^\fg
  \]
  Since $\dim V_{\varpi_2} = 6$, we have that $\twedge^{2k} V_{\varpi_2}=0$ for $k>3$. We have the following
  isomorphisms of $\fg$-modules
  \[
    \twedge^2 V_{\varpi_2} \cong \twedge^4 V_{\varpi_2} \cong V_{\varpi_1+\varpi_3},\qquad
    \twedge^6 V_{\varpi_2} \cong \Cee.
  \]
  Therefore, $\exp\circ\mu = 1+\mu + \tfrac12\mu^2 + \tfrac16\mu^3$.

  $\bullet$ For $(\fg,V) = (\fo_n,V_{\varpi_1})$, $W=V_{\varpi_1}$, where $n>6$ we have that
  \begin{equation}\label{eq:isoPhiOn}
    \left(\cF_{\poly}^{2k}\otimes (W\otimes W^\ast)\right)^\fg
    \cong  \left(\twedge^{2k}V_{\varpi_1} \otimes (V_{\varpi_2}\oplus V_{2\varpi_1}\oplus \Cee)\right)^\fg    
  \end{equation}
  Using~\eqref{eq:OnE4V1}, we get that~\eqref{eq:isoPhiOn} vanishes for $k>1$. Therefore, $\exp\circ\mu = \id + \mu$.

  Finally, we see that since $\mu$ is nilpotent in all case, $\exp\circ \mu$ has the maximal rank.
\end{proof}

\ssbegin{Lemma}\label{lem:LMoment}
Under assumptions of Theorem~\ref{thm:main}, the bivector $\pi_B - r_{\cV}$ is Poisson and the $G$-action on
$(\cV, \pi_B -  r_{\cV})$ is Poisson. Moreover,
\begin{enumerate}
\item
  for $B$ symmetric, the $G$-action  admits a~moment map in sense of Lu
  $L = (L_{+}, L_{-}): \cV \to G^\ast$, where $\Phi = L_{+}L_{-}^{-1}$;
  
  \item 
  for $B$ symplectic, 
  the pre-image under $\Phi$ of the domain of the Gauss decomposition carries a Hamiltonian $\fg$-actions and a moment
  map in the sense of Lu $L = (L_{+}, L_{-}):  \cV \to G^\ast$, where $\Phi = L_{+}L_{-}^{-1}$,
  \item the Poisson structures $\pi_B$ and $-r_{\cV}$ on~$\cV$ are compatible.
\end{enumerate}
\end{Lemma}
\begin{proof}
  Consider a~twist of~$\pi' = \pi_B+r_{\dyn}(\mu)$ by the antisymmetrized $r$-matrix~$r$. 
  It follows from Lemma~\ref{lem:PhiFactorL}
  that $(\cV,\pi_B+r_{\dyn}(\mu)-r_{\cV})$ is a Poisson space
  and the Gauss decomposition of~$G$ defines a~$G^\ast$-valued moment map~$L =
 (L_{+}, L_{-})$ by factorising $\Phi = L_{+}L_{-}^{-1}$. 
  In particular, we get that
  \[
    [\![\pi_B + r_{\dyn}(\mu) - r_{\cV}, \pi_B + r_{\dyn}(\mu) - r_{\cV}]\!] = 0.
  \]
  By Lemma~\ref{lem:Rdyn},  $r_{\dyn}(\mu) = 0$ and this implies
  $$
  [\![\pi_B - r_{\cV}, \pi_B - r_{\cV}]\!] = 0
  $$
proving the first claim. Note that by~\cite[Theorem~3.12]{Zwicknagl2009} 
 $[\![r_{\cV},r_{\cV}]\!]=0$ for $(\fg,V)$ listed in Lemma~\ref{lem:LIST}, and this  shows that  $[\![\pi_B,r_{\cV}]\!] = 0$.

\end{proof}

\ssbegin{Remark}
We can assume that $r$ is given by the anti-symmetrization of any $r$-matrices from~\cite{BelavinDrinfeld}.
Then, for $(L_{+},L_{-})\in G^\ast$ we have that $L_{+}L_{-}^{-1} \in G$, see~\cite{SemenovTianShansky1985}.
If $\Phi$ belongs to the image of such a factorisation,  then in the
same manner we get a $G^\ast$-valued moment map.
\end{Remark}

\section{Examples for $\fg=\fsl_2$ }
In this section, we consider examples of even and odd  modules for~$\fg=\fsl_2$.
Let $e$, $h$, $f$ be the standard basis of~$\fsl_2$ and
$B_\fg$ be the nondegenerate symmetric invariant bilinear form on~$\fsl_2$ is given by
\begin{equation}\label{eq:sl2NIS}
  B_\fg(e,f) = 1, \quad B_\fg(h, h)  = 2.
\end{equation}
Therefore, the corresponding Cartan 3-tensor takes the form
\[
  \phi =  \frac12 e \wedge h \wedge f.
\]

Recall the standard antisymmetrised $r$-matrix for~$\fsl_2$ 
\[
  r = \frac 12 e \wedge f 
\]
and the corresponding cobracket 
\[
  \delta(h) = 0,\qquad \delta(e) = - \frac12 h\wedge e,\qquad
  \delta(f) = - \frac12 h\wedge f.
\]
The dual Poisson Lie group~$SL(2,\Cee)^\ast$ is defined by
\[  
  SL(2,\Cee)^\ast = \left\{ L_+ =
    \begin{pmatrix} a & b \\ 0 & a^{-1} \end{pmatrix},\quad
    L_- \begin{pmatrix} a^{-1} & 0 \\ c & a \end{pmatrix} \mid \text{$a,b,c\in\Cee$ such that $a\neq 0$}\right\}.
\]
Note that $L_+L_{-}^{-1}\in SL(2,\Cee)$.
The Lie algebra $\fsl(2)^\ast$ can be realised as pairs of matrices
\[
  \fsl(2)^\ast = \left\{
    \begin{pmatrix} a/2 & b \\ 0 & -a/2 \end{pmatrix},
    \begin{pmatrix} -a/2 & 0 \\ -c & a/2 \end{pmatrix}
    \mid a,b,c\in\Cee
  \right\},
\]
with a basis $h^\vee = \frac12(h,-h)$, $e^\vee=(e,0)$, $f^\vee = (0,-f)$.

\subsection{Motivating example: $\twedge \fsl_2$}
Let us consider an odd $\fsl_2$-module $V = V_{2\varpi_1}$ isomorphic to the adjoint representation.
We view the algebra~$\cF=\twedge\fsl_2$ as an~algebra of function on a $0|3$-dimensional
linear supermanifold~$\cV$.
Let $\xi_2$, $\xi_0$, $\xi_{-2}$ be the basis in~$\twedge^1\fsl_2$ corresponding to the standard
basis $e$, $h$, $f$ of~$\fsl_2$.
The Lie algebra homomorphism $\fsl_2 \to \fX(\cV)$ is given by
\begin{align*}
  e_{\cV} = {}&\xi_0\frac{\del}{\del \xi_{-2}} - 2\xi_2\frac{\del}{\del\xi_0},
  & h_{\cV} ={}& 2\xi_2\frac{\del}{\del \xi_2} - 2\xi_{-2}\frac{\del}{\del\xi_{-2}},
  & f_{\cV} ={}& -\xi_0\frac{\del}{\del\xi_2} + 2\xi_{-2}\frac{\del}{\del\xi_0}.
\end{align*}

The form~$B$ on~$V$ is induced by~\eqref{eq:sl2NIS}.
The bivector
\[
  \pi_B = 2\frac{\del}{\del\xi_2}\wedge\frac{\del}{\del\xi_{-2}} + 2\frac{\del}{\del\xi_0}\wedge\frac{\del}{\del\xi_0} \in \mathbb{E}^2\fX(\cV)
\]
corresponds to the Poisson bracket on~$\cF$ defined by
\[
  \{\xi,\eta\}_B = 2B(\xi,\eta)\qquad
  \text{for $\xi,\eta\in\twedge^1\fsl_2\subset\cF$}
\]
The quadratic moment map is of the form
\[
  \mu = \frac12\left(\xi_2\xi_{-2}\otimes h^\ast - \xi_0\xi_{-2}\otimes e^\ast -\xi_2\xi_0\otimes f^\ast \right).
\]

It is easy to check that
\[
  \phi_{\cV} = \frac{1}{2} e_{\cV}\wedge h_{\cV}\wedge f_{\cV} = 0.
\]
Hence, we can view $(\cV,\pi_B)$ as a quasi-Poisson supermanifold.
Since $\twedge^4 V = 0$, we get that
\[
  r_{\dyn}(\mu) \in (\mathbb{E}^{2,4}\fX_{\poly}(\cV))^\fg
  \oplus (\mathbb{E}^{2,8}\fX_{\poly}(\cV))^\fg\oplus\ldots
  \cong (\twedge^4V\otimes S^2(V))^\fg \oplus\ldots
  = 0.
\]
Therefore, the correction term $r_{\dyn}(\mu)$ given by the dynamical $r$-matrix vanishes.

The direct computations show that for the quasi-Poisson
structure $\pi_{\Phi} = \pi_B + r_{\dyn}(\mu) = \pi_B$ on~$\cV$
the $SL(2,\Cee)$-valued moment map~$\Phi$ is given by
\[
  \Phi = \exp(\mu) = \begin{pmatrix}
    1 + \frac12 \xi_2\xi_{-2} & - \frac12 \xi_0\xi_{-2} \\
    - \frac12 \xi_2\xi_0 & 1 - \frac12 \xi_2 \xi_{-2}
  \end{pmatrix}.
\]

The anti-symmetrized $r$-matrix $r\in\twedge^2\fsl_2$  defines a Poisson bracket on~$\cF$ by
\[
  \{a, b \} = ({\rm mult} \circ r)(a \otimes b))
  \qquad\text{for $a,b \in\cF$}.
\]
The corresponding bivector is given by
\[
  r_{\cV} =
  \xi_0\xi_{-2}\frac{\del}{\del\xi_0}\wedge\frac{\del}{\del\xi_{-2}} - 2\xi_2\xi_{-2}\frac{\del}{\del\xi_0}\wedge\frac{\del}{\del\xi_0} +  \xi_2 \xi_0\frac{\del}{\del\xi_2}\wedge\frac{\del}{\del\xi_0}.
\]
The bivector $ \pi_B- r_{\cV}$ corresponds to the twist of~$\pi_B$ by~$r$, see~\S\ref{sec:twist}.
The Gauss decomposition of~$SL(2,\Cee)$ leads to the following factorisation of~$\Phi$:
\[
  \Phi =
    \begin{pmatrix}
      1 & 0 \\ \frac12 \xi_0\xi_2 & 1
    \end{pmatrix}
    \begin{pmatrix}
      1 + \frac12 \xi_2\xi_{-2} & 0 \\ 0 & 1 - \frac12 \xi_2\xi_{-2}
    \end{pmatrix}
    \begin{pmatrix}
      1 & -\frac12 \xi_0\xi_{-2} \\ 0 & 1
    \end{pmatrix}.    
\]
Therefore, we get a $SL(2,\Cee)^\ast$-valued moment map $L = (L_{+}, L_{-})$, where
\[
  L_{+} = \begin{pmatrix}
    1 + \frac14 \xi_2\xi_{-2} & - \frac12 \xi_0\xi_{-2} \\
    0 & 1 - \frac14 \xi_2\xi_{-2}
  \end{pmatrix},\qquad
  L_{-} = \begin{pmatrix}
    1 - \frac14 \xi_2\xi_{-2} & 0 \\
    \frac12 \xi_2\xi_0 &  1+ \frac14\xi_2\xi_{-2}
  \end{pmatrix}.
\]

Finally, we note that there is a change of variables
\[
  \xi_2 \mapsto \xi_2,\qquad
  \xi_0 \mapsto \xi_0 + \xi_2\xi_0\xi_{-2},\qquad
  \xi_{-2} \mapsto \xi_{-2}
\]
mapping the  bivector $\pi_B$ to $\pi_B - r_{\cV}$.
This change of variables is not $\fsl_2$-invariant. Therefore, the map~$\mu$ is not a~moment map
for the Poisson structure~$\pi_B - r_{\cV}$. 
Observe that this change of variables is related to a
non-standard filtration on the Clifford algebra of~$\fsl_2$, see Remark~3.6 in~\cite{cubicDiracUqSL2}.
Quantisation of~$\twedge\fsl_2$ and of the $G^\ast$-valued moment map~$L$ with respect to the Poisson structure given by $\pi_B-r_{\cV}$ were considered in~\cite{qSL2Diff,cubicDiracUqSL2}.

\subsection{The case of $(\fg,V) = (\fsl_2,V_{\varpi_1})$}

The module $V=V_{\varpi_1}$ is even and admits a~nondegenerate symplectic bilinear form given by
\[
  B(v_1, v_{-1}) = 1,\quad B(v_{-1}, v_{1}) = -1,\quad
  B(v_1, v_1) = B(v_{-1}, v_{-1})  = 0,
\]
where $v_1$, $v_{-1}$ is the standard basis of~$V$ and $v_1^\ast = -v_{-1}$, $v_{-1}^\ast = v_1$ is the corresponding
dual basis with respect to~$B$. The action of~$\fg$ on~$\cV$ is given by the following vector fields
\begin{align*}
  e_{\cV} = {}& v_{1}\frac{\del}{\del v_{-1}},
  & h_{\cV} = {}& v_1\frac{\del}{\del v_1} - v_{-1}\frac{\del}{\del v_{-1}},
  & f_{\cV} = {}& v_{-1}\frac{\del}{\del v_{1}}.
\end{align*}
The Poisson structure induced by~$B$ is given by the bivector
\[
  \pi_B = 2 \frac{\del}{\del v_1}\wedge\frac{\del}{\del v_{-1}}.
\]
The moment map is of the form
\[
  \mu =
  \frac14 v_{1}^2\otimes e^\ast  - \frac12 v_1v_{-1}\otimes h^\ast - \frac14 v_{-1}^2\otimes f^\ast
  \in \cF\otimes\fg^\ast.
\]

For dimensional reasons, we have that
\[
  \phi_{\cV} =  \frac12 e_{\cV}\wedge h_{\cV} \wedge f_{\cV}  =0.
\]
Thus, any invariant Poisson structure on~$\cV$ is also quasi-Poisson.
It is easy to check that $r_{\dyn}(\mu) = 0$.

Recall that $e^\ast = f$, $h^\ast = \tfrac12 h$, $f^\ast = e$ with respect to~\eqref{eq:sl2NIS}. We get
\[
  \mu = \frac14 \begin{pmatrix}
    - v_1v_{-1} & - v_{-1}^2 \\
    v_{1}^2 &  v_1v_{-1}
  \end{pmatrix}\in \cF^2_{\poly}\otimes\fg.
\]
A direct calculation shows that $\mu^2=0$.
Therefore, the expression 
\[
  \Phi = \exp\circ\mu =  E + \mu = E + \frac14 \begin{pmatrix}
    - v_1v_{-1} & - v_{-1}^2 \\
    v_1^2 & v_1v_{-1}
  \end{pmatrix}
\]
defines an $SL(2,\Cee)$-valued moment map for $(\cV,\pi_B)$.

The bivector~$r_{\cV}$ corresponding to the Poisson bracket~\eqref{eq:brR} on~$\cV$ given by the antisymmetrized
$r$-matrix takes the form
\[
  r_{\cV} = - \frac12 v_1 v_{-1} \frac{\del}{\del v_1}\wedge \frac{\del}{\del v_{-1}}.
\]
The Gauss decomposition of~$\Phi$ is given by
\[
  \Phi =
  \begin{pmatrix} 1 & \frac{-v_{-1}^2}{4+ v_1v_{-1}}\\ 0 & 1 \end{pmatrix}
  \begin{pmatrix} \frac{4}{4 + v_1v_{-1}} & 0 \\ 0 &  \frac14(1+4v_1v_{-1})  \end{pmatrix}
  \begin{pmatrix} 1 & 0 \\\frac{v_1^2}{4 + v_1v_{-1}} & 1  \end{pmatrix}.
\]
Assuming that $1+\frac14v_1v_{-1}>0$, we get that
\[
  L_{+} = \begin{pmatrix}
    \frac{1}{\sqrt{1+\frac14 v_1v_{-1}}} & - \frac{v_{-1}^2}{4\sqrt{1+\frac14 v_1v_{-1}}} \\
    0 & \sqrt{1+\frac14 v_1v_{-1}}
  \end{pmatrix},\quad
  L_{-} = \begin{pmatrix}
    \sqrt{1+\frac14 v_1v_{-1}} & 0 \\
    - \frac{v_1^2}{4\sqrt{1+\frac14v_1v_{-1}}}& \frac{1}{\sqrt{1+\frac14 v_1v_{-1}}}
  \end{pmatrix}    
\]
are such that $\Phi = L_{+} L_{-}^{-1}$. Therefore, we have a $SL(2,\Cee)^\ast$-valued moment map given by~$L=(L_{+},L_{-})$.

\section{Bivectors on $V \oplus V^\ast$}
The following question is motivated by the definition of $q$-deformed Clifford algebras in~\cite{cubicDiracUqSL2}.
Let $V$ be $\fg$-module such that $\Hom_\fg(\twedge^3V,S^3(V))=0$.
Then, the same is true for its dual $V^\ast$.
Thus, the actions of $G$ on~$(\cV,-r_{\cV})$ and $(\cV^\ast,-r_{\cV^\ast})$ are Poisson.
Let $W = V\oplus V^\ast$ and define on it a natural invariant symmetric or symplectic bilinear form~$B$ by the duality between~$V$ and~$V^\ast$.

\ssbegin{Problem}\label{pr:fusion}
Construct a~Poisson bivector $\pi_{\cW}$ on the supermanifold $\cW$ corresponding to~$W$ such that
\begin{enumerate}
\item $\pi_{\cW}$ is quadratic;

\item the action of $G$ on $(\cW,\pi_{\cW})$ is Poisson;

\item the bivector $\pi_{\cW}$ is compatible with $\pi_B$;

\item the restriction of $\pi_{\cW}$ to $\cV$ is $-r_{\cV}$ and the restriction to~$\cV^\ast$ is $-r_{\cV^\ast}$.
\end{enumerate}
\end{Problem}

Note that the bivectors $ (-r_{\cV} - r_{\cV^\ast})$ and $(-r_{\cW})$ do not satisfy the conditions stated above since
\[
  \cL_{x_{\cW}}(-r_{\cV} - r_{\cV^\ast})
  = - \delta(x)_{\cV} - \delta(x)_{\cV^\ast}  \neq - \delta(x)_{\cW},
  \hskip 0.3cm [[ -r_{\cW}, -r_{\cW} ]] \neq  0.
\]

We answer the question stated above in the case of $\fg=\fsl_n$ and $V=V_{\varpi_1}$.
Recall that the element
\[
  \psi = \frac12 \sum e_a^1 \wedge f_a^2 \in \twedge^2(\fg\oplus\fg),
\]
where $e_a^1$ is a~basis of~$\fg$ embedded in~$\fg\oplus\fg$ as the first copy
and $f_a^2$ is the dual basis (with respect to~$B_\fg$) of~$\fg$ embedded in~$\fg\oplus\fg$ as the second copy,
plays an~important role in the fusion of the quasi-Poisson manifolds, see~\cite[\S5]{AlekseevKosmannSchwarzbachMeinrenken2002}.
For the rest of this Section, we use $v_i$ for a~basis of~$V$ and $v_i^\ast$ for the dual basis of~$V^\ast$
regardless of their parity. Set $\del_i := \frac{\del}{\del v_i}$ and $\del_i^\ast := \frac{\del}{\del v_i^\ast}$.

\ssbegin{Theorem}\label{thm:SLfusion}
Let $(\fg,V) = (\fsl_n,V_{\varpi_1})$. Let $v_i$ be the standard basis of~$V$ and $v_i^\ast$ be the
corresponding dual basis of~$V^\ast$. Then, the bivector
\begin{align*}
  \pi_\cW ={}& - r_{\cW} - \psi_{\cW} - \left(-(-1)^s\frac{1}{2} - \frac{1}{2\dim V}\right)\sum_{i,j} v_i\del_i \wedge v_j^\ast \del^\ast_j\\
  \intertext{
  or equivalently}
  {}={}&  - r_{\cW} - \psi_{\cW}
  + \frac{1}{2\dim V}\sum_{i,j} v_i\del_i \wedge v_j^\ast \del^\ast_j
         + \frac{1}{2}\sum_{i,j} v_j^\ast v_i\del_i \wedge  \del^\ast_j,
\end{align*}
satisfies
\[
  [[\pi_{\cW},\pi_{\cW}]] = 0,\qquad
  [[\pi_B,\pi_{\cW}]] = 0,\qquad
  L_{x_{\cW}} \pi_{\cW} = - \delta(x)_{\cW}.
\]
Here $s=1$ if we view $W$ as an odd module, and $s=0$ if we view it as an even module,
\end{Theorem}

\begin{proof}
  Set
  \[
    \pi_{cor} = \sum_{i,j} v_i\del_i \wedge v^\ast_j\del^\ast_j.
  \]
  Since $\psi_{\cW}$ and $\pi_{cor}$ are $G$-invariant, we have that  
  \begin{align*}
    \cL_{x_{\cW}}( \pi_{\cW})  = -\cL_{x_{\cW}}( r_{\cW}) = - \delta(x)_{\cW}.
  \end{align*}  
  This implies that
  \[
    [\![ r_{\cW}, \psi_{\cW} ]\!] = 0,\quad
    [\![ r_{\cW}, r_{\cW} ]\!] = -\phi_{\cW},\quad
    [\![\psi_{\cW},\psi_{\cW}]\!] =  \phi_{\cW},\quad
    [\![r_{\cW}, \pi_{cor} ]\!] = 0,\quad
    [\![\psi_{\cW}, \pi_{cor} ]\!] = 0,\quad
  \]
  which leads to
  \[
    [\![ \pi_{\cW} , \pi_{\cW} ]\!] = \phi_{\cW} - \phi_{\cW}  = 0.
  \]
  This proves that the action of~$G$ on~$(\cW,\pi_\cW)$ is Poisson.

  For $[\![\pi_B,\pi_{\cW}]\!]=0$, note that $V$ and $V^\ast$ are also $\fgl_n$-module.
  Let $B_\fg$ be the trace form on~$\fgl_n$, then the restriction of~$B_\fg$ on~$\fsl_n$ is nondegenerate.
  The matrix unit $E_{i,j}\in\fgl_n$ is dual to $E_{j,i}$ with respect to~$B_\fg$. 
Let $\psi^\fgl \in \twedge^2(\fgl_n\oplus \fgl_n)$, then we have that
\[
  \psi^\fgl = \psi + \frac{1}{2\dim V} z^1\wedge z^2,
\]
where $z=\sum_i E_{i,i}$ is the identity matrix in~$\fgl_n$.
The factor $\dim V$ in the denominator comes from the fact that $B_\fg(z,z) = \dim V$.
We have that
\[
  (E_{i,j})_{\cV} = v_i\del_j,\qquad
  (E_{i,j})_{\cV^\ast} = - v^\ast_j \del^\ast_i,\qquad
  (z)_{\cV} = \sum_i v_i\del_i,\qquad
  (z)_{\cV^\ast} = -\sum_j v^\ast_j\del^\ast_j.  
\]
Hence $\pi_{cor} =  - (z^1\wedge z^2)_{\cW}$.
Therefore, we have that
\[
  \psi_{\cW} = \psi^\fgl_{\cW} - \frac{1}{2\dim V} (z^1\wedge z^2)_{\cW}
  = \psi^\fgl_{\cW} + \frac{1}{2\dim V} \pi_{cor},
\]
\[
  \psi^\fgl_{\cW} = \frac12 \sum_{i,j} (E_{i,j})_{\cV}\wedge(E_{j,i})_{\cV^\ast}
    = -\frac12 \sum v_i\del_j\wedge v^\ast_i\del^\ast_j.
\]
We have
\begin{align*}
  [\![\pi_B, \pi_{cor}]\!] ={}
  & [\![2\del_i\wedge \del^\ast_i,v_a\del_a\wedge v_b^\ast\del^\ast_b]\!] \\
  {}={}& 2 ([\![\del_i\wedge \del^\ast_i,v_a\del_a]\!]\wedge v_b^\ast\del^\ast_b
         - v_a\del_a\wedge [\![\del_i\wedge \del^\ast_i,v_b^\ast \del^\ast_b]\!])\\
  {}={}&  2([\![\del_i,v_a\del_a]\!]\wedge\del^\ast_i\wedge v_b^\ast\del^\ast_b
         - v_a\del_a\wedge \del_i\wedge [\![\del^\ast_i,v_b^\ast \del^\ast_b]\!])\\
  {}={}&  2(\del_i\wedge\del^\ast_i\wedge v_b^\ast\del^\ast_b
         - v_a\del_a\wedge \del_i\wedge \del^\ast_i)\\
  {}={}&  2(v_b^\ast\del_i\wedge\del^\ast_i\wedge \del^\ast_b
         - v_a\del_a\wedge \del_i\wedge \del^\ast_i).
\end{align*}
On the other hand, we get
\begin{align*}
  [\![\pi_B,\psi^{\fgl}_{\cW}]\!] = {}
  &  [\![2\del_i\wedge \del^\ast_i, \frac12(E_{ab})_{\cV}\wedge (E_{ba})_{\cV^\ast} ]\!]\\
  {}={}& [\![\del_i\wedge \del^\ast_i, v_a\del_b\wedge (-v^\ast_a\del^\ast_b)]\!]\\
  {}={}&- ([\![\del_i\wedge \del^\ast_i, v_a\del_b]\!]\wedge v^\ast_a\del^\ast_b
         - v_a\del_b \wedge [\![\del_i\wedge \del^\ast_i, v^\ast_a\del^\ast_b]\!])\\
  {}={}&-( [\![\del_i,v_a\del_b]\!]\wedge \del^\ast_i \wedge v^\ast_a\del^\ast_b
         - v_a\del_b \wedge \del_i \wedge [\![\del^\ast_i, v^\ast_a\del^\ast_b]\!])\\
  {}={}&-(\del_b \wedge \del^\ast_i \wedge v^\ast_i \del^\ast_b
         - v_i\del_b\wedge \del_i \wedge \del^\ast_b).\\
  \intertext{If $\cW$ is an odd module, then}
  [\![\pi_B,\psi^{\fgl}_{\cW}]\!] = {}
  & -( v^\ast_i \del_b\wedge \del^\ast_b \wedge \del^\ast_i
    - v_i\del_i\wedge \del_b \wedge \del^\ast_b).
\end{align*}
Finally, we get that
\begin{align*}
  [\![\pi_B,\pi_{\cW}]\!] ={}
  & [\![\pi_B,-\psi_{\cW}-c\pi_{cor} ]\!]
  = - [\![\pi_B,\psi^\fgl_{\cW} + \frac{1}{2\dim V}\pi_{cor} + c\pi_{cor}]\!]\\
  {}={}&-(-1 + \frac{1}{\dim V} + 2c)(
         v_b^\ast\del_i\wedge\del^\ast_i\wedge \del^\ast_b
         - v_a\del_a\wedge \del_i\wedge \del^\ast_i)
\end{align*}
Therefore, $c = \frac{\dim V-1}{2\dim V} = \frac12 - \frac{1}{2\dim V}$ for odd modules.

If $\cW$ is an even module, then
\[
  [\![\pi_B,\psi^{\fgl}_{\cW}]\!] = {}
   v^\ast_i \del_b\wedge \del^\ast_b \wedge \del^\ast_i
    - v_i\del_i\wedge \del_b \wedge \del^\ast_b.
\]
Finally, we get that
\begin{align*}
  [\![\pi_B,\pi_{\cW}]\!] ={}
  & [\![\pi_B,-\psi_{\cW}-c\pi_{cor} ]\!]
  = - [\![\pi_B,\psi^\fgl_{\cW} + \frac{1}{2\dim V}\pi_{cor} + c\pi_{cor}]\!]\\
  {}={}&-(1 + \frac{1}{\dim V} + 2c)(
         v_b^\ast\del_i\wedge\del^\ast_i\wedge \del^\ast_b
         - v_a\del_a\wedge \del_i\wedge \del^\ast_i)
\end{align*}
Therefore, $c = -\frac{\dim V+1}{2\dim V}= - \frac12 - \frac{1}{2\dim V}$ for even modules.
\end{proof}

\ssbegin{Remark}
The Poisson structure $(\cW,\pi_\cW)$ described in Theorem~\ref{thm:SLfusion} is related to
the quantisation of the algebra of differential forms on the Grassmannians, see~\cite{HKdR} and~\cite{NicholsGrass}.
\end{Remark}

\subsection{Examples}
In this subsection we consider a~couple of examples of different moment maps related to the Poisson structure
on~$\cW$ described in Theorem~\ref{thm:SLfusion}. In what follows we assume that $W = V\oplus V^\ast$ is
an odd module.
\subsubsection{$(\fsl_2,V_{\varpi_1}\oplus V_{\varpi_1}^*)$}
Let $V$ be the tautological $\fsl_2$-module and consider $\cW = V \oplus V^\ast$ as an odd module by fixing
an~invariant symmetric bilinear~$B$ form~$\cW$. Then the corresponding linear supermanifold~$\cW$ is a
Hamiltonian Poisson $SL(2,\Cee)$-supermanifold with the Poisson bivector given by
\[
  \pi_B = 2\del_1\wedge\del_1^\ast + 2\del_2\wedge \del_2^\ast
\]
and the moment map defined by
\[
  \mu  = \frac14(v_1v_1^\ast - v_2v_2^\ast)\otimes h
  + \frac12 v_2v_1^\ast \otimes e
  + \frac12 v_1v_2^\ast \otimes f
\]
here $v_i\in V$  and $v_i^\ast\in V^\ast$ are odd generators of~$\cF$ such that
\begin{gather*}
  e v_1 = 0,\qquad h v_1 =  v_1,\qquad f v_1 = v_2,\qquad
  e v_2 = v_1,\qquad h v_2 =  -v_1,\qquad f v_2 = 0,\\
  e v_1^\ast = -v_2^\ast,\qquad h v_1^\ast =  -v_1^\ast,\qquad f v_1^\ast = 0,\qquad
  e v_2^\ast = 0,\qquad h v_2^\ast =  v_2^\ast,\qquad f v_2^\ast = v_1^\ast.
\end{gather*}
The corresponding to~$\mu$ dynamical $r$-matrix is given by
\[
  r_{\dyn}(\mu)  = - \frac18 v_1v_2v_1^\ast v_2^\ast \del_1\del_1^\ast
   - \frac18 v_1v_2v_1^\ast v_2^\ast \del_2\del_2^\ast. 
 \]
Then, $(\cW, \pi_B + r_{\dyn}(\mu))$ is a~quasi-Poisson supermanifold which admits the $SL(2,\Cee)$-valued moment
map given by
\[
  \Phi  = \exp\circ\mu = \begin{pmatrix}
    1 + \frac14 v_1v_1^\ast - \frac14 v_2v_2^\ast + \frac{3}{16} v_1v_2v_1^\ast v_2^\ast & \frac12 v_2v_1^\ast \\
    \frac12 v_1v_2^\ast  & 1 - \frac14 v_1v_1^\ast + \frac14 v_2v_2^\ast + \frac{3}{16} v_1v_2v_1^\ast v_2^\ast
  \end{pmatrix}.
\]
The twist given by the standard $r$-matrix is defined by
\[
  r_{\cW} = \frac12( v_1v_2\del_1\wedge\del_2 - v_1v_1^\ast\del_2\wedge \del_2^\ast
  +v_2v_2^\ast\del_1\wedge\del_1^\ast - v_1^\ast v_2^\ast \del_1^\ast\wedge\del_2^\ast).
\]
The Poisson bivector from Theorem~\ref{thm:SLfusion} is given by
\begin{align*}
  \pi_{\cW}  
  &{}= v_1v_1^\ast \del_2\wedge \del_2^\ast + \frac12(- v_1v_2 \del_1\wedge \del_2 - v_1v_2^\ast\del_1\wedge\del_2^\ast
  - v_2v_1^\ast \del_2\wedge\del_1^\ast + v_1^\ast v_2^\ast \del_1^\ast\wedge\del_2^\ast).
\end{align*}
Let us consider the new quasi-Poisson structure on~$\cW$ which is given by the bivector
\[
  \pi_B + \pi_{\cW} + r_{\cW}.
\]
The direct computations show the following surprising fact: the $SL(2,\Cee)$-valued moment map for $(\cW,
\pi_B + \pi_{\cW} + r_{\cW})$ is the same as for $(\cW, \pi_B + r_{\dyn}(\mu))$!

Since the moment map $\Phi$ is given by the exponent of~$\mu$, we see that the quasi-Poisson supermanifold
$(\cW, \pi_B + \pi_{\cW} + r_{\cW})$ is the exponent of the Hamiltonian Poisson supermanifold~$\cW$ with the bivector defined by
\[
  \pi_B + \pi_{\cW} + r_{\cW} + r_{\dyn}(\mu)
\]
and the moment map defined by~$\mu$.

It is easy to see that both $(\cW,\pi_B)$ and $(\cW, \pi_B + \pi_{\cW} + r_{\cW} + r_{\dyn}(\mu))$ are symplectic
supermanifolds. The corresponding symplectic structures are given by
\[
  \omega_B = 2 \dd v_1\wedge \dd v_1^\ast + 2 \dd v_2\wedge \dd v_2^\ast
\]
and
\[
  \omega = \dd v_1\wedge \dd v_1^\ast (2 - \frac12 v_1v_2^\ast)  + \frac12 \dd v_1\wedge \dd v_2^\ast   v_2v_1^\ast
  + \dd v_2\wedge \dd v_3 \frac12 v_1v_2^\ast  + \dd v_2\wedge \dd v_2^\ast (2 - \frac12 v_1v_1^\ast).
\]
Note that here we consider the differential forms on~$\cW$ as a~right $\cF$-module.

The difference of symplectic forms is the differential of an invariant basic 1-form
\[
  \omega - \omega_B  = \frac12 \dd\left( \dd v_1 v_2v_1^\ast v_2^\ast - \dd v_2 v_1v_1^\ast v_2^\ast\right).
\]
Hence, for any $x\in\fsl_2$ we get that
\[
  \dd \langle\mu,x \rangle = \iota_{X_\cW}\omega_B
  = \iota_{X_{\cW}}\left(\omega - \frac12 \dd\left( \dd v_1 v_2v_1^\ast v_2^\ast - \dd v_2 v_1v_1^\ast v_2^\ast\right)\right)
  = \iota_{X_{\cW}}\omega.
\]
Therefore, $\mu$ serves as a~moment map for both Poisson structures.

\subsubsection{The case of $(\fg,V) = (\fsl(3), V_{\varpi_1})$}
Already in this case, $\Phi = \exp\circ \mu$ doesn't define a~moment map for the quasi-Poisson structure on~$\cW$
given by the bivector
\[
  \pi_B + \pi_{\cW} + r_{\cW}.
\]
The $SL(3,\Cee)$-valued moment for this bivector is given by $\exp\circ\nu$, where
\allowdisplaybreaks
\begin{align*}
  \nu = {}
  & \left(\frac{1}{3} v_1 v^*_1 - \frac{1}{6} v_2 v^*_2 - \frac{1}{6} v_3 v^*_3
    +\frac{1}{24} v_1 v_2 v^*_1 v^*_2 + \frac{1}{24} v_1 v_3 v^*_1 v^*_3 - \frac{1}{12} v_2 v_3 v^*_2  v^*_3\right)\otimes H_1\\
  {}&{} + \left(\frac{1}{6} v_1 v^*_1 + \frac{1}{6} v_2 v^*_2 - \frac{1}{3} v_3 v^*_3
      +\frac{1}{12} v_1 v_2 v^*_1 v^*_2-\frac{1}{24} v_1 v_3 v^*_1 v^*_3-\frac{1}{24} v_2 v_3 v^*_2  v^*_3\right)\otimes H_2\\
  {}&{} + \left(\frac{1}{2} v_2 v^*_1 +\frac{1}{8} v_2 v_3 v^*_1 v^*_3\right)\otimes E_{1,2}
      - \left(\frac{1}{2} v_3 v^*_2-\frac{1}{8} v_1 v_3 v^*_1 v^*_2\right)\otimes E_{2,3}\\
  {}&{} - \left(\frac{1}{2} v_3 v^*_1+\frac{1}{8} v_2 v_3 v^*_1 v^*_2\right)\otimes E_{1,3}
      + \left(\frac{1}{2} v_1 v^*_2 +\frac{1}{8} v_1 v_3 v^*_2 v^*_3\right)\otimes E_{2,1}\\
  {}&{}  + \left(-\frac{1}{2} v_2 v^*_3  -\frac{1}{8} v_1 v_2 v^*_1 v^*_3\right)\otimes E_{3,2}
      + \left(\frac{1}{2} v_1 v^*_3-\frac{1}{8} v_1 v_2 v^*_2 v^*_3\right)\otimes E_{3,1},
\end{align*}
where $E_{i,j}$ are matrix units and $H_1 = E_{1,1} - E_{2,2}$, $H_2 = E_{2,2} - E_{3,3}$.

Recall that the moment map for~$\pi_B$ is given by
\begin{align*}
  \mu = {}
  & \left(\frac{1}{3} v_1 v^*_1 - \frac{1}{6}v_2 v^*_2 - \frac{1}{6}v_3 v^*_3\right)\otimes H_1
 + \left(\frac{1}{6} v_1 v^*_1 + \frac{1}{6}v_2 v^*_2-\frac{1}{3}v_3 v^*_3\right)\otimes H_2\\
  {}&{}  + \frac{1}{2}v_2 v^*_1\otimes E_{1,2}  -\frac{1}{2}v_3 v^*_2\otimes E_{2,3}
      - \frac{1}{2} v_3 v^*_1 \otimes E_{1,3} \\
{}&{} + \frac{1}{2} v_1 v^*_2\otimes E_{2,1}
    - \frac{1}{2} v_2 v^*_3\otimes E_{3,2}
    + \frac{1}{2} v_1v^*_3\otimes E_{3,1}.
\end{align*}
We note that the quadratic terms in~$\mu$ and~$\nu$ are the same.

\providecommand{\bysame}{\leavevmode\hbox to3em{\hrulefill}\thinspace}
\providecommand{\MR}{\relax\ifhmode\unskip\space\fi MR }
\providecommand{\MRhref}[2]{%
  \href{http://www.ams.org/mathscinet-getitem?mr=#1}{#2}
}
\providecommand{\href}[2]{#2}


\end{document}